\documentclass[12pt]{article}
\usepackage{graphicx,psfrag,epsfig, color}
\usepackage{amssymb,amsmath,amscd,amsthm}
\usepackage{graphicx,psfrag,epsfig}

\usepackage{graphicx}
\usepackage[active]{srcltx}

\newtheorem{theorem}{Theorem}[section]

\newtheorem{lemma}[theorem]{Lemma}

\date{}

\begin{document}

\date{}
\title{On the behavior of diffusion processes with traps}
\author{
 M. Freidlin\footnote{Dept of Mathematics, University of Maryland,
College Park, MD 20742, mif@math.umd.edu}, L.
Koralov\footnote{Dept of Mathematics, University of Maryland,
College Park, MD 20742, koralov@math.umd.edu}, A. Wentzell\footnote{Dept of Mathematics, Tulane University, New Orleans, LA  70118, wentzell@math.tulane.edu}
} \maketitle

\begin{abstract}
We consider processes that coincide with a given diffusion process outside a finite collection of domains. In each of the domains, there is, additionally, a large drift directed towards the interior of the domain. We describe the limiting behavior of the processes as the magnitude of the drift tends to infinity, and thus the domains become trapping, with the time to exit the domains being exponentially large. In particular, in exponential time scales, metastable distributions between the trapping regions are considered.
\end{abstract}

{2010 Mathematics Subject Classification Numbers: 60F10, 35J25, 47D07, 60J60.}

{ Keywords: Metastable Distributions, Large Deviations, Exit Problem, Non-standard Boundary Problem.  }

\section{Introduction} \label{intro}

Let $v$ be an infinitely smooth vector field on the $d$-dimensional torus $\mathbb{T}^d$. (General manifolds, whether compact or not, can also
be considered, but we'll stick with the torus for the sake of simplicity of notations later on.) Consider the process $X^{x, \varepsilon}_t$ defined via
\begin{equation} \label{mainpro}
d X^{x, \varepsilon}_t = \frac{1}{\varepsilon} v ( X^{x, \varepsilon}_t) d t + d W_t,~~~X^{x, \varepsilon}_0 = x,
\end{equation}
where $\varepsilon$ is a small parameter and $W$ is a $d$-dimensional Brownian motion. This process is an order-one perturbation of a strong deterministic flow, or, equivalently, a time-changed small random perturbation of the deterministic flow $\dot{x}(t) = v(x(t))$. We'll be interested in the large-time behavior of the process when $\varepsilon \downarrow 0$.

More precisely, assume at first that the collection of asymptotically stable limit sets of the unperturbed flow consists of $n$ equilibrium points $O_1$,...,$O_n$. Let $D_1$,...,$D_n$ denote the sets that are attracted to $O_1$,...,$O_n$, respectively. Intuitively, when $\varepsilon$ is small and $t$ is large, $X^{x, \varepsilon}_t$ is located near one of these points with overwhelming probability. The transitions between small neighborhoods of the equilibriums are governed by the matrix
\begin{equation} \label{vij}
V_{ij} = \frac{1}{2} \inf \{  \int_0^T | \dot{\varphi}_s - v(\varphi_s)|^2 ds,~~\varphi_0 = O_i, \varphi_T = O_j\},~~i,j \in \{1,...,n\},
\end{equation}
where the infimum is taken over all $T \geq 0$ and all absolutely continuous functions   $\varphi: [0,T] \rightarrow \mathbb{T}^d$. Loosely speaking, for each $i = 1,...,n$, $x \in  D_i$, and   $\lambda > 0$  (except a finite subset of $\lambda$'s, as discussed below)  there is an index $k = k(x,\lambda) \in \{1,...,n\}$ such that
\begin{equation} \label{mesta}
{\rm dist}(X^{x, \varepsilon}_{\exp(\lambda/\varepsilon)}, O_k) \rightarrow 0~~~{\rm in}~{\rm probability}~{\rm as}~\varepsilon \downarrow 0.
\end{equation}
The equilibrium $O_k$ is called the metastable state for the process $X^{x, \varepsilon}_t$ corresponding to the initial point $x$ and time scale $\exp(\lambda/\varepsilon)$. The state $O_k$ can be determined by comparing $\lambda$ with certain linear expressions involving the numbers $V_{ij}$ (see Chapter 6 of \cite{FW}). Typically, there is a finite set $\Lambda$ of values of $\lambda$ where transition from one metastable state to another happens, i.e., the notion of a metastable state is defined for $\lambda \in (0,\infty)\setminus \Lambda$.

On the other hand, for certain geometries of the unperturbed flow, it may happen that $V_{ij_1} = V_{ij_2}$ for some $j_1 \neq j_2$ or, more generally, the sums of two distinct collections of $V_{ij}$'s may be equal. The analysis of the asymptotics of $X^{x, \varepsilon}_t$  is then more intricate. The notion of rough symmetry, of which the simplest examples are due to geometric symmetries of the flow, was introduced and to some extend analyzed in \cite{FF}. In the presence of rough symmetry, metastable states may need to be replaced by metastable distributions between the asymptotically stable attractors.

\begin{figure}[htbp]
%
%
\centerline{\includegraphics[height=2.5in, width= 3.3in,angle=0]{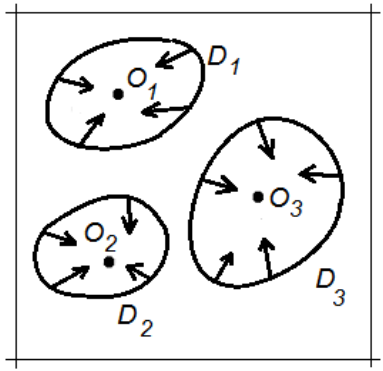}}
  \vskip -20pt
  \caption{Torus with multiple trapping regions.}
    \label{domains}
\end{figure}

One of the motivations of the current paper is to analyze an interesting situation where $V_{ij}$ do not depend on $j$ and, consequently, $X^{x, \varepsilon}_t$ is distributed between the neighborhoods of several equilibriums at exponential time scales. Now this phenomenon is due not to geometric symmetries but to vanishing of the vector field in certain regions of the state space.  Namely, let $D_1,...,D_n$ be open connected domains with infinitely smooth boundaries $\partial D_k$, $k =1,...,n$, on the $d$-dimensional torus $\mathbb{T}^d$. The closures $\overline{D}_k$ are assumed to be disjoint. Let $v$ be a vector field on $\mathbb{T}^d$ that is equal to zero on $ \mathbb{T}^d \setminus \bigcup_{k=1}^n \overline{D}_k$. It is assumed to be infinitely smooth in the sense that there is an infinitely smooth field on $ \mathbb{T}^d$ that
agrees with $v$ on  $\bigcup_{k=1}^n \overline{D}_k$. Assume, for the moment, that all the points of $D_k$ are attracted to an equilibrium $O_k \in D_k$, $k =1,...,n$ (see Figure~\ref{domains}).
\vskip -10pt
  \vskip 10pt

The quantities $V_{ij}$ are easily seen not to depend on $j$ since
\[
\inf \{ \int_0^T | \dot{\varphi}_s - v(\varphi_s)|^2 ds,~~\varphi_0 \in \partial D_i, \varphi_T \in \partial D_j\} =
\]
\[
\inf \{  \int_0^T | \dot{\varphi}_s|^2 ds,~~\varphi_0 \in \partial D_i, \varphi_T \in \partial D_j\} = 0,
\]
where the infimum is taken over all $T > 0$ and $\varphi \in C^1([0,T], \mathbb{T}^d)$. There are two  issues involved in  understanding the transitions of the process $X^{x, \varepsilon}_t$ between the neighborhoods of different equilibriums. The first issue is to describe how the process that starts near $O_k$ exits the domain $D_k$. Fortunately, this questions has been well-studied. One needs to look at the quasi-potential
\begin{equation} \label{quap}
V_k(x) = \frac{1}{2} \inf \{  \int_0^T | \dot{\varphi}_s - v(\varphi_s)|^2 ds,~~\varphi_0 = O_k, \varphi_T = x\},~~x \in \partial D_k,
\end{equation}
where the infimum is taken over all $T \geq 0$ and all $\varphi \in C^1([0,T], \overline{D}_k)$. If the minimum of $V_k(x), x \in \partial D_k$,
is achieved in a single point $x_k$, then the process starting near $O_k$ exits $D_k$ in a small neighborhood of $x_k$. The time to exit is of order $\exp(V_k(x_k)/\varepsilon)$. Even if the minimum is not achieved in a single point, there are cases when the exit from the domain is well understood (e.g., in the simplest example when $v$ is sphericaly symmetric and $D_k$ is a ball around $O_k$, also see~\cite{Day} and references there).
We'll simply assume that for each compact $K \subset D_k$, the exit time (appropriately re-scaled) and the exit location have limiting distributions that do not depend on the starting point within $K$, as in the case of a single minimum for the quasi-potential and in the symmetric case mentioned above.

The second issue concerns the transition between the domains $D_k$, i.e., the behavior of the process  $X^{x, \varepsilon}_t$ on $ \mathbb{T}^d \setminus \bigcup_{k=1}^n {D}_k$. In order to get a meaningful limiting object, we introduce a new process $Y^{x, \varepsilon}_t$ by running the clock only when $X^{x, \varepsilon}_t \in \mathbb{T}^d \setminus \bigcup_{k=1}^n {D}_k$. While this process obviously coincides with the Brownian motion away from boundary $\bigcup_{k=1}^n \partial {D}_k$, the regions $D_k$ still play an important role by capturing the process when it reaches the boundary and then re-distributing it on $\partial D_k$. The description of the limiting behavior $Y^{x, \varepsilon}_t$ is the main result of this paper.

In Section~\ref{dttp} we describe the limiting process.  It belongs to a peculiar class of processes that, it seems, has not been discussed in the literature previously. In Section~\ref{seconv} we prove the convergence of $Y^{x, \varepsilon}_t$ to the limit. In Section~\ref{appli} we
describe the asymptotics of the original process $X^{x, \varepsilon}_t$  at exponential time scales and make several additional remarks.

\section{Description of the limiting process} \label{dttp}

In this section we define the family of processes $X^x_t$, which later will be proved to be the limiting processes for $Y^{x, \varepsilon}_t$  as $\varepsilon \downarrow 0$.
Let $D_1,...,D_n \subset \mathbb{T}^d$ be open connected domains with infinitely smooth boundaries $\partial D_k$, $k =1,...,n$. The closures $\overline{D}_k$ are assumed to be disjoint. Let $ U = \mathbb{T}^d \setminus \bigcup_{k=1}^n \overline{D}_k$. The closure of this domain will be denoted by $ \overline{U}$. Let $U'$ be the metric space obtained from $\overline{U}$ by identifying all points of $\partial D_k$, making every $\partial D_k$, $k =1,...,n$, into one point $d_k$.

The family of processes $X^x_t$, $x \in U'$,  will be defined in terms of its generator. Since we expect $X^x_t$ to coincide with a Brownian motion inside $U$, the generator coincides with $\frac{1}{2} \Delta$ on a certain class of functions. The domain, however, should be restricted by certain boundary conditions to account for non-trivial behavior of $X^x_t$ on the boundary of $U$. We'll use the Hille-Yosida theorem stated here in the form that is convenient for considering closures of linear operators.
\begin{theorem}
Let $K$ be a compact space, $C(K)$ be the space of continuous functions on it.  Suppose that a linear operator $A$ on $C(K)$ has the following properties:

(a) The domain $\mathcal{D}(A)$ is dense in $C(K)$;

(b) The identity $ \mathbf{1}$ belongs to $\mathcal{D}(A)$ and $A \mathbf{1} = 0$;

(c) The maximum principle: If $S$ is the set of points where a function $f \in \mathcal{D}(A)$ reaches its maximum, then $A f (x) \leq 0$ for at
least one point $x \in S$.

(d) For a dense set $\Psi \subseteq C(K)$ for every $\psi \in \Psi$ and every $\lambda > 0$ there exists a solution $f \in \mathcal{D}(A)$ of the equation $\lambda f - A f = \psi$.

Then the operator $A$ is closable and its closure $ \overline{A}$ is the infinitesimal generator of a unique semi-group of positivity-preserving operators $T_t$, $t \geq 0$, on $C(K)$ with $T_t \mathbf{1} = \mathbf{1}$, $||T_t|| \leq 1$.
\end{theorem}

Suppose that we are given positive finite measures $\nu_1,...,\nu_n$ concentrated on $\partial D_1$,...,$\partial D_n$, respectively.
The Hille-Yosida theorem will be applied to the space $K = U'$. Let us define the linear operator $A$ in $C(U')$. First we define its domain: it consists of all functions $f \in C(U')$ that satisfy the following conditions:

(1) $f$  is twice continuously differentiable in $U$;

(2) The limits of all the first and second order derivatives of $f$ exist at all the points of the boundary $\partial U = \bigcup_{k=1}^n \partial {D}_k$;

(3) There are constants $g_1,...,g_n$ such that
\[
\lim_{y \in U, {\rm dist}(y, \partial D_k) \downarrow 0} \Delta f(y) = g_k,~~~k =1,...,n;
\]

(4)  For each $k =1,...,n$,
\begin{equation} \label{intco}
\int_{\partial D_k} \langle \nabla f(x) , n(x) \rangle \nu_k(d x) =0,
\end{equation}
where  $n(x)$ is the unit exterior (with respect to $U$) normal at $x \in \partial D_k$.

For $f \in \mathcal{D}(A)$ and $x \in U'$, we define
\[
A f =  \begin{cases} \frac{1}{2} \Delta f(x), & {\rm if }~ x \in U, \\ \frac{1}{2} g_k, & {\rm if }~ x = d_k,~~~k=1,...,n. \end{cases}
\]
Let us check that the conditions of the Hille-Yosida theorem are satisfied.

(a) Consider the set $G$ of functions $g$ that are infinitely smooth and have the following property: for each $k =1,...,n$ there is
a set $V_k$ open in $U'$ such that $\partial D_k \subset V_k$ and $g$ is constant on $V_k$. It is clear that $G \subset \mathcal{D}(A)$ and $G$ is dense in $C(U')$.

(b) Clearly $\mathbf{1} \in \mathcal{D}(A)$ and $A \mathbf{1} = 0$.

(c) If $f$ has a maximum at $x \in U$, it is clear that $ \Delta f (x) \leq 0$. Now suppose that $f$ has a maximum at $d_k$. We can view $f$ as an element of  $C^2(\overline{U})$ that is constant on each component of the boundary, in particular on $\partial D_k$. Note that $\langle \nabla f(x) , n(x) \rangle$ is identically zero on
$\partial D_k$, since otherwise it would be negative at some points due to (\ref{intco}), which would contradict the fact that $f$ reaches its maximum on $\partial D_k$. Then the second derivative of $f$ in the direction of $n$ is non-positive at all points $x \in \partial D_k$. Since $f$ is constant on $\partial D_k$, its second derivative in any direction tangential to the boundary is equal to zero. Therefore, $\Delta f(x) \leq 0$ for $x \in \partial D_k$, i.e., $A f (d_k) \leq 0$, as required.

(d) Let $\Psi$ be the set of functions $\psi \in C(U')$ that have limits of all the first order derivatives as $y \in U, y \rightarrow x$, at all points
$x \in \partial U$. It is clear that $\Psi$ is dense in $C(U')$. Let $\widetilde{f} \in C^2 (\overline{U})$ be the solution of the equation
$\lambda \widetilde{f} - \frac{1}{2} \Delta \widetilde{f} = \psi$ in $U$, $\widetilde{f} = 0$ on $\partial U$. Let $h_k \in C^2(\overline{U})$ be the solution of the equation
\[
\lambda {h}_k(x) - \frac{1}{2} \Delta {h}_k(x) = 0,~~~~x \in U,
\]
\[
h_k(x) = 1,~~x \in \partial D_k;~~~~~h_k(x) = 0,~~x \in \partial U \setminus \partial D_k.
\]
Let us look for the solution $f \in \mathcal{D}(A)$ of $\lambda f -  A f = \psi$ in the form $f = \widetilde{f} + \sum_{k =1}^n c_k h_k$. We get $n$ linear equations for $c_1,...,c_n$. The solution is unique because of the maximum principle. Therefore, the determinant of the system is non-zero, and the solution exists for all the right hand sides.

Let $ \overline{A}$ be the closure of $A$. Let $T_t$, $t \geq 0$, be the corresponding semi-group on $C(U')$, whose existence is guaranteed by the
Hille-Yosida theorem. By the Riesz-Markov-Kakutani representation theorem, for $x \in U'$ there is a measure $P(t,x,dy)$ on $(U', \mathcal{B}(U'))$ such that
\[
(T_t f)(x) = \int_{U'} f(y)P(t,x,dy),~~~f \in C(U').
\]
It is a probability measure since $T_t \mathbf{1} = \mathbf{1}$. Moreover, it can be easily verified that $P(t,x,B)$ is a Markov transition function.  Let $X^x_t$, $x \in U'$, be the corresponding Markov family. In order to show that a modification with continuous trajectories exists, it is enough to check that $\lim_{t \downarrow 0} P(t,x,B)/t = 0$ for each closed set $B$ that doesn't contain $x$ (Theorem I.5 of \cite{Mandl}, see also \cite{Dyn}). Let  $f \in \mathcal{D}(A)$  be a non-negative function that is equal to one on $B$ and whose support doesn't contain $x$. Then
\[
\lim_{t \downarrow 0} \frac{P(t,x,B)}{t} \leq \lim_{t \downarrow 0} \frac{(T_t f)(x) - f(x)}{t} = Af(x) = 0,
\]
as required. Thus $X^x_t$ can be assumed to have continuous trajectories.

\section{Convergence of the trace of the process} \label{seconv}

In this section we prove the convergence of $Y^{x, \varepsilon}_t$, obtained from $X^{x, \varepsilon}_t$ by running the clock only when the process is in $\overline{U}$, to the limiting process $X^x_t$.

Let $v$ be a vector field that is smooth in $ \bigcup_{k=1}^n \overline{D}_k$ (i.e., it admits a smooth continuation from $ \bigcup_{k=1}^n \overline{D}_k$ to the whole space) and is equal to zero outside $ \bigcup_{k=1}^n \overline{D}_k$.

Let
\[
a(x) = \langle v(x), n(x)\rangle,~~x \in \partial U,
\]
where $n$ is the unit exterior (with respect to $U$) normal  to the boundary. We'll assume that $a(x) > 0$ for all $x \in \partial U$.
Recall that the process $X^{x, \varepsilon}_t$ is defined via
\[
d X^{x, \varepsilon}_t = \frac{1}{\varepsilon} v ( X^{x, \varepsilon}_t) d t + d W_t,~~~X^{x, \varepsilon}_0 = x.
\]
For $B \subset \mathbb{R}$, let $  \tau^{x, \varepsilon}(B) =  \inf \{t\geq 0: X^{x,\varepsilon}_t \in B \}$.
Let $\mu_k^{x, \varepsilon} $ be the measure on $\partial D_k$ induced  by $X^{x, \varepsilon}_{\tau^{x, \varepsilon}(\partial D_k)}$. We'll  assume that
there are measures $\mu_k$, $k=1,...,n$, such that for each compact set $K \subset D_k$ and each continuous function $\varphi$ on $\partial D_k$ we have
\begin{equation} \label{convqq}
\lim_{\varepsilon \downarrow 0} \int_{\partial D_k} \varphi d \mu_k^{x, \varepsilon} = \int_{\partial D_k} \varphi d \mu_k
\end{equation}
uniformly in $x \in K$.
Define
\[
s(t) = \inf(s: \lambda(u: u \leq s, X^{x, \varepsilon}_u \in \overline{U}) > t)),
\]
where $\lambda$ is the Lebesgue measure on the real line, and
\begin{equation} \label{stc}
Y^{x, \varepsilon}_t = X^{x, \varepsilon}_{s(t)}.
\end{equation}
Thus $Y^{x, \varepsilon}_t $ is a right-continuous process with values in $\overline{U}$, which also can be viewed as a right-continuous $U'$-valued process. It can be obtained from $X^{x, \varepsilon}_t $ by running the clock only when $X^{x, \varepsilon}_t $ is in $\overline{U}$.   Define the measures $\nu_k$ via
\begin{equation} \label{meas1}
\nu_k(dx) = (2 a(x))^{-1} \mu_k(dx),~~~x \in \partial D_k.
\end{equation}

Let $X^x_t$ be the Markov family of continuous $U'$-valued processes defined above, corresponding to the measures $\nu_k$, $k =1,...,n$. The main result of this section is the following.
\begin{theorem} \label{mt1}
For each $x \in U'$, the measures on $C([0,\infty), U')$ induced by the processes $Y^{x, \varepsilon}_t $ converge weakly, as $\varepsilon \downarrow 0$, to the measure induced by $X^x_t$.
\end{theorem}
 The key ingredient in the proof of this theorem is the following lemma.
\begin{lemma}
Suppose that $f \in \mathcal{D}(A)$. Then
\begin{equation} \label{maeq}
\lim_{\varepsilon \downarrow 0} \mathrm{E} \left(f(Y^{x, \varepsilon}_t) - f(x) -  \frac{1}{2} \int_0^t \Delta f (Y^{x, \varepsilon}_u) du \right) = 0
\end{equation}
for each $t \geq 0$, uniformly in  $x \in \overline{U}$.
\end{lemma}
\proof
For the sake of notational simplicity, we'll assume that there is just one domain where the vector field $v$ is non-zero. This does not lead to
any loss of generality as the proof in the case of multiple domains is similar. We'll denote the domain by $D$ and will drop the subscript $k$ from the notations everywhere. For example, (\ref{intco}) now takes the form
\begin{equation} \label{ortt}
\int_{\partial D} {\langle \nabla f (x), n(x) \rangle} d \nu(x) = 0.
\end{equation}

Let $S_r = \{x \in \overline{U}: {\rm dist}(x, \partial D) = r\}$ for $r \geq 0$,  $S_r = \{x \in D: {\rm dist}(x, \partial D) = -r\}$ for $r < 0$. These are smooth surfaces if $r$ is sufficiently small. Let $\Gamma_r = \{x \in \overline{U}: {\rm dist}(x, \partial D) \leq r \}$ for $r \geq 0$.

Let $\sigma^{x, \varepsilon}_0 = 0$,  $\tau^{x, \varepsilon}_1 = \tau^{x, \varepsilon}(\partial D)$,  $\sigma^{x, \varepsilon}_n = \inf(t \geq \tau^{x, \varepsilon}_n: X^{x, \varepsilon}_t \in S_{\sqrt{\varepsilon}})$, $n \geq 1$, while $\tau^{x, \varepsilon}_n = \inf(t \geq \sigma^{x, \varepsilon}_{n-1}: X^{x, \varepsilon}_t \in \partial D)$, $n \geq 2$. Then
\[
\mathrm{E} \left(f(Y^{x, \varepsilon}_t) - f(x) -  \frac{1}{2}  \int_0^t \Delta f (Y^{x, \varepsilon}_u) du \right) =
\]
\begin{equation} \label{vv1}
\mathrm{E} \sum_{n = 1}^\infty  \left(f(X^{x, \varepsilon}_{\tau^{x, \varepsilon}_n \wedge s(t)} ) - f(X^{x, \varepsilon}_{\sigma^{x, \varepsilon}_{n-1} \wedge s(t)}) - \frac{1}{2}  \int_{\sigma^{x, \varepsilon}_{n-1} \wedge s(t)}^{\tau^{x, \varepsilon}_n \wedge s(t)} \Delta f (X^{x, \varepsilon}_u) du \right) +
\end{equation}
\[
\mathrm{E} \sum_{n = 1}^\infty  \left(f(X^{x, \varepsilon}_{\sigma^{x, \varepsilon}_n \wedge s(t)} ) - f(X^{x, \varepsilon}_{\tau^{x, \varepsilon}_{n} \wedge s(t)}) - \frac{1}{2}  \int_{\tau^{x, \varepsilon}_{n} \wedge s(t)}^{\sigma^{x, \varepsilon}_n \wedge s(t)} \Delta f (X^{x, \varepsilon}_u) du \right),
\]
where we put $\Delta f \equiv 0 $ on $D$. The first expectation on the right hand side is equal to zero since $X^{x, \varepsilon}_t$ is pure Brownian motion on $\overline{U}$. Next, we prove several statements primarily concerning the stopping times $ \sigma^{x, \varepsilon}_n$ and the behavior of the process
near $\partial D$.
\\

1) There is $c= c(t) > 0$ such that
\begin{equation} \label{jj1}
 \mathrm{P}(\sigma^{x, \varepsilon}_n  \leq s(t)) \leq \exp({-{ c n \sqrt{\varepsilon}}}),~~~x \in \overline{U},~ n \geq 2.
\end{equation}
Indeed, since $\partial D$ is smooth, there is $r > 0$ such that the ball of radius $r$ tangent to $\partial D$ at $x$ lies entirely in $\overline{U}$. Let $\eta^\varepsilon$ be the time it takes a Brownian motion starting inside a ball of radius $r$ at a distance $\sqrt{\varepsilon}$ from the boundary to reach the boundary. It is easy to see that there is $c = c(t) > 0$ such that
\[
 \mathrm{P}(\eta^\varepsilon  \leq t) \leq \exp({-{ c  \sqrt{\varepsilon}}}).
\]
Therefore, if $\eta^\varepsilon_k$, $k \geq 1$, is a sequence of such independent random variables, then
\begin{equation} \label{jj1s}
 \mathrm{P}(\eta^\varepsilon_1 + ... + \eta^\varepsilon_n  \leq t) \leq \exp({-{c n \sqrt{\varepsilon}}}),~~~ n \geq 1.
\end{equation}
By our construction, $\mathrm{P}(\tau^{x, \varepsilon}_n - \sigma^{x, \varepsilon}_{n-1} > z| X^{x, \varepsilon}_{\sigma^{x, \varepsilon}_{n-1}}) \geq
\mathrm{P}(\eta^\varepsilon > z)$ for each $n \geq 2$ and $z>0$. Therefore, estimate (\ref{jj1}), with $\tau^{x, \varepsilon}_{n+1}$ instead of $ \sigma^{x, \varepsilon}_n$, follows from (\ref{jj1s}) and the strong Markov property. Thus, the original formula (\ref{jj1}) also holds, with a different constant $c$.
\\

2) Consider an auxiliary one-dimensional process $Z^z_t$ that satisfies
\[
d Z^z_t = d B_t - {\mathrm{v}}  \chi_{(-\infty, 0)}(Z^z_t) dt,~~~Z^z_0 = z,
\]
where $ \mathrm{v} > 0$ and $B_t$ is a one-dimensional Brownian motion. Also consider its perturbation defined via
\begin{equation} \label{dife}
d \widetilde{Z}^{z}_t = d B_t - {\mathrm{v}}  \chi_{(-\infty, 0)}(\widetilde{Z}^{z}_t) dt + {A^z_t} dt ,~~~\widetilde{Z}^{z}_0 = z,
\end{equation}
where $A^z_t$ is a continuous adapted process satisfying $|A^z_t| \leq \mathrm{v}/2$.  It is easy to see that for each $\eta > 0$ there
is $z_0 < 0$  such that
\begin{equation} \label{ny1}
\mathrm{P}(\sup_{t \geq 0} \widetilde{Z}^{z_0}_t \geq 0) \leq \eta,
\end{equation}
while for each $z_0 < 0$ there is $t_0$ such that for $z \in [z_0, 1]$,
\begin{equation} \label{yy1}
\mathrm{P}({\widetilde{Z}}^{z}_t~{\rm reaches}~\{z_0\}  \cup \{1\}~{\rm before}~{\rm time}~t_0) \geq 1-\eta.
\end{equation}
 A direct calculation shows that
\[
\mathrm{P}(\sup_{t \geq 0} {Z}^0_t \geq 1) = (1+ 2  \mathrm{v})^{-1}.
\]
Fix such $z_0$ that (\ref{ny1}) holds and
\begin{equation} \label{cl1}
\mathrm{P}( {Z}^0_t~{\rm reaches}~ 1 ~{\rm before}~{\rm reaching}~ z_0) \geq (1+ 2  \mathrm{v})^{-1} - \eta.
\end{equation}
By the Girsanov formula, whose use is justified by (\ref{yy1}) with ${\widetilde{Z}}^{z}_t$ replaced by ${{Z}}^{0}_t$, there is $\varkappa >0$ such that
\begin{equation} \label{fvv}
\mathrm{P}( {\widetilde{Z}}^{0}_t~{\rm reaches}~ 1 ~{\rm before}~{\rm reaching}~ z_0) \in
[(1+ 2  \mathrm{v})^{-1} - 3\eta, (1+ 2  \mathrm{v})^{-1} + 3\eta],
\end{equation}
provided that $|A^0_t| \leq \varkappa$.
Below we'll encounter a related $\varepsilon$-dependent process. Namely, suppose that $(\widetilde{Z}^{z,\varepsilon}_t, {\widehat{Z}}^{z,\varepsilon}_t) \in \mathbb{R} \times \mathbb{R}^{d-1}$ satisfy
\begin{equation} \label{mij}
d \widetilde{Z}^{z,\varepsilon}_t = d \widetilde{B}_t - \frac{\mathrm{v}}{\varepsilon}  \chi_{(-\infty, 0)}(\widetilde{Z}^{z,\varepsilon}_t) dt + \frac{\widetilde{A}^z_t}{\varepsilon} dt ,~~~\widetilde{Z}^{z,\varepsilon}_0 = \widetilde{z},
\end{equation}
\begin{equation} \label{mijxx}
d {\widehat{Z}}^{z,\varepsilon}_t = \sigma(\widetilde{Z}^{z,\varepsilon}_t, {\widehat{Z}}^{z,\varepsilon}_t) d \widehat{B}_t + \frac{{\widehat{A}}^z_t}{\varepsilon} dt ,~~~{\widehat{Z}}^{z,\varepsilon}_0 = {\widehat{z}},
\end{equation}
where $\widetilde{B}_t$ is a one-dimensional Brownian motion, $\widehat{B}_t$ is $d$-dimensional Brownian motion, possibly correlated with $\widetilde{B}_t$, $\sigma$ is a $(d-1) \times d$ matrix, and $z = (\widetilde{z}, \widehat{z}) \in \mathbb{R}^d$. We'll assume that there is $C > 0$ such that
\begin{equation} \label{ez}
|\widetilde{A}^z_t| \leq \mathrm{v}/2,~~~~|\sigma|, | \widehat{A}^z_t| \leq C.
\end{equation}
For $A \subset \mathbb{R}$, let $ \widetilde{\tau}^{z, \varepsilon}(A)  =  \inf \{t\geq 0: \widetilde{Z}^{z,\varepsilon}_t \in A \}$.
We claim that there are $\varepsilon_0 >0$ and $L > 0$ such that
\begin{equation} \label{us2}
\mathrm{E}\left( \lambda(t: {\widetilde{Z}}^{0,\varepsilon}_t \in [0, \varepsilon], t \leq \widetilde{\tau}^{0, \varepsilon}(\{-\sqrt{\varepsilon}\} \cup \{\varepsilon \}))\right) \leq L \varepsilon^{2},
\end{equation}
provided that $\varepsilon \leq \varepsilon_0$.
%
If, additionally,
\begin{equation} \label{mij2}
|\widetilde{A}^0_t| \leq \varkappa(\varepsilon)~~{\rm whenever}~|\widetilde{Z}^{0,\varepsilon}_t| + |\widehat{Z}^{0,\varepsilon}_t| \leq \sqrt{\varepsilon},~~{\rm with}~\varkappa(\varepsilon) \downarrow 0,~{\rm as}~\varepsilon \downarrow 0,
\end{equation}
then
\begin{equation} \label{us1}
\lim_{\varepsilon \downarrow 0} \mathrm{P}\left( \widetilde{\tau}^{0, \varepsilon}(\{\varepsilon\}) < \widetilde{\tau}^{0, \varepsilon}(\{-\sqrt{\varepsilon}\})\right) =
(1+ 2  \mathrm{v})^{-1}.
\end{equation}
Indeed, let $\eta \in ( 0,1)$. Find $z_0 < 0$ such that (\ref{cl1}) holds and
\begin{equation} \label{iip}
\mathrm{P}(\widetilde{\tau}^{z_0 \varepsilon, \varepsilon}(\{0\}) < \infty) \leq \eta,
\end{equation}
which is possible by (\ref{ny1}) and the scaling-invariance of the Brownian motion.
By (\ref{yy1}), we can find $t_0$ such that for $z \in [z_0, 1]$,
\begin{equation} \label{kk1}
\mathrm{P}\left(\widetilde{\tau}^{z \varepsilon, \varepsilon}(\{z_0 \varepsilon\} \cup \{\varepsilon\}) \leq t_0 \varepsilon^2\right)
\geq 1-\eta.
\end{equation}
The combination of these two inequalities and the strong  Markov property (whose use is allowed since a time-shift of the process $\widetilde{A}^z_t$ is also bounded by $ \mathrm{v}/2$)  imply (\ref{us2}).
By (\ref{mij2}), (\ref{kk1}), (\ref{mij}), and (\ref{mijxx}),
\[
\mathrm{P}( \sup_{t \leq \widetilde{\tau}^{0, \varepsilon}(\{z_0 \varepsilon \} \cup \{\varepsilon\})} |\widetilde{A}^0_t| \leq \varkappa(\varepsilon)) \geq
\mathrm{P}(\sup_{t \leq \widetilde{\tau}^{0, \varepsilon}(\{z_0 \varepsilon \} \cup \{\varepsilon\})} (|\widetilde{Z}^{0,\varepsilon}_t| + |\widehat{Z}^{0,\varepsilon}_t|) \leq \sqrt{\varepsilon} ) \geq 1 - 2 \eta
\]
for all sufficiently small $\varepsilon$. Therefore, by (\ref{fvv}), which can be applied after re-scaling,
\[
\mathrm{P}\left( \widetilde{\tau}^{0, \varepsilon}(\{\varepsilon\}) < \widetilde{\tau}^{0, \varepsilon}(\{z_0 \varepsilon\})\right) \in
[(1+ 2  \mathrm{v})^{-1} - 5\eta, (1+ 2  \mathrm{v})^{-1} + 5\eta].
\]
Together with (\ref{iip}), this implies that
\[
\mathrm{P}\left( \widetilde{\tau}^{0, \varepsilon}(\{\varepsilon\}) < \widetilde{\tau}^{0, \varepsilon}(\{-\sqrt{\varepsilon}\})\right) \in
[(1+ 2  \mathrm{v})^{-1} - 6\eta, (1+ 2  \mathrm{v})^{-1} + 6\eta].
\]
Since $\eta$ was arbitrary, this implies (\ref{us1}).

From the presence of a strong drift to the left in (\ref{mij}), it easily follows that there is $L > 0$ such that
\begin{equation} \label{q2}
\sup_{z \in [-\sqrt{\varepsilon}, \varepsilon] } \mathrm{E} \widetilde{\tau}^{z, \varepsilon}(\{-\sqrt{\varepsilon}\} \cup \{\varepsilon\}) \leq L \varepsilon^{\frac{3}{2}}.
\end{equation}
From (\ref{q2}), (\ref{mij}),  and (\ref{mijxx}) it follows that
\begin{equation} \label{iui}
\lim_{\varepsilon \downarrow 0} \mathrm{P}(\sup_{t \leq \widetilde{\tau}^{0, \varepsilon}(\{-\sqrt{\varepsilon} \} \cup \{\varepsilon\})} (|\widetilde{Z}^{0,\varepsilon}_t| + |\widehat{Z}^{0,\varepsilon}_t|) > \varepsilon^{\frac{1}{3}} ) = 0.
\end{equation}
\\

3) 
Let us show that for each sufficiently small $\delta > 0$  there are  $\varepsilon_0 >0$ and $L > 0$ such that
\begin{equation} \label{us2xg}
\mathrm{E}\left( \lambda(t: X^{x,\varepsilon}_t \in \overline{U}, t \leq \min( \tau^{x,\varepsilon}(S_{-\delta}), \tau^{x,\varepsilon}(S_\varepsilon)))\right) \leq L \varepsilon^{2},~~~x \in \partial D,
\end{equation}
provided that $\varepsilon \leq \varepsilon_0$, and that
\begin{equation} \label{us1xg}
\lim_{\varepsilon \downarrow 0} \mathrm{P}\left( \tau^{x,\varepsilon}(S_\varepsilon) < \tau^{x,\varepsilon}(S_{-\delta})\right) =
(1+ 2  a(x))^{-1}~~~{\rm uniformly}~{\rm in}~x \in \partial D.
\end{equation}
First observe that
\begin{equation} \label{interm}
\lim_{\varepsilon \downarrow 0} \varepsilon^{-2} \mathrm{P}( \tau^{x,\varepsilon}(S_\varepsilon) < \tau^{x,\varepsilon}(S_{-\delta})) = 0~~~{\rm uniformly}~{\rm in}~x \in S_{-\sqrt{\varepsilon}}~,
\end{equation}
due to the presence of the strong drift inside $D$. Also, there is $L > 0$ such that
\begin{equation} \label{brmo}
\mathrm{E} \min( \tau^{x,\varepsilon}(\partial D), \tau^{x,\varepsilon}(S_\varepsilon)) \leq L \varepsilon^{2},~~~x \in \Gamma_\varepsilon,
\end{equation}
for all sufficiently small $\varepsilon$, since the process is pure Brownian motion on $U$.


Let us describe a change of coordinates in a neighborhood of a point $x \in \partial D$. Let $V^r_\varepsilon = [-\sqrt{\varepsilon}, \varepsilon] \times B_r \subset \mathbb{R}^d$, where $B_r \subset \mathbb{R}^{d-1}$ is the closed ball of radius $r$ centered at the origin. Let $m_x$ be an isometric mapping of $B_r$ to the $(d-1)$-dimensional ball of radius $r$ centered at $x$ in the tangent plane to $\partial D$ at $x$. For $y \in B_r$,
we take the straight line passing through $m_x(y)$ and the point on $\partial D$ closest to $m_x(y)$ (this line is perpendicular to $\partial D$ if $m_x(y) \notin \partial D$; we define it as the perpendicular if $m_x(y) \in \partial D$).  For $z \in [-\sqrt{\varepsilon}, \varepsilon]$ and $y \in B_r$, define $\varphi_x(z, y)$ as the point on the perpendicular that belongs to $S_z$.
If $r$ and $\varepsilon$ are sufficiently small, then $\varphi_x$ is a diffeomorphism from ${V}^r_\varepsilon$ to a domain ${U}^r_\varepsilon(x)$ for each~$x$.

For $z \in {V}^r_\varepsilon$, let $\overline{X}^{z, \varepsilon}_t = \varphi^{-1}({X}^{\varphi(z), \varepsilon}_t)$ be the process written in the new coordinates (stopped when it reaches the boundary of ${V}^r_\varepsilon$). It satisfies
\[
d \overline{X}^{z, \varepsilon}_t = \frac{1}{\varepsilon} \overline{v} ( \overline{X}^{z,\varepsilon}_t)  \chi_{[-\sqrt{\varepsilon},0] \times \mathbb{R}^{d-1}}( \overline{X}^{z,\varepsilon}_t) d t + \overline{\beta} ( \overline{X}^{ z,\varepsilon}_t) d t + \overline{\sigma} (\overline{X}^{z, \varepsilon}_t) d \overline{W}_t,~~~\overline{X}^{z,\varepsilon}_0 = z.
\]
The coefficients $\overline{v}, \overline{\beta}, \overline{\sigma}$ are bounded in $C^1(V^r_\varepsilon)$ (the change of coordinates, and therefore the coefficients, depend on $x$, but the bound is uniform in $x$) and satisfy: $\overline{v}_1(0) = -a(x)$, $\overline{\sigma}(0)$ is an orthogonal matrix. Let
\[
\alpha(z) = (\sigma_{11}^2(z)+...+\sigma_{1d}^2(z))^{-\frac{1}{2}}.
\]
Note that this is a smooth function such that $\alpha(0) = 1$. The process
\[
d {\widetilde{X}}^{z, \varepsilon}_t = \frac{1}{\varepsilon} \alpha^2 \overline{v}   \chi_{[-\sqrt{\varepsilon},0] \times \mathbb{R}^{d-1}}( {\widetilde{X}}^{z,\varepsilon}_t) d t + \alpha^2 \overline{\beta} ( {\widetilde{X}}^{z, \varepsilon}_t) d t + \alpha \overline{\sigma} ({\widetilde{X}}^{z, \varepsilon}_t) d \overline{W}_t,~~~{\widetilde{X}}^{z, \varepsilon}_0 = z,
\]
is different from $\overline{X}^{z, \varepsilon}_t$ by a random change of time. The coefficients of $ {\widetilde{X}}^{z, \varepsilon}_t$ can be extended from $V^r_\varepsilon = [-\sqrt{\varepsilon}, \varepsilon] \times B_r \subset \mathbb{R}^d$ to $[-\sqrt{\varepsilon},\varepsilon] \times \mathbb{R}^{d-1}$ as continuous functions by requiring that they don't vary in the radial direction outside $B_r$. This way the process ${\widetilde{X}}^{z, \varepsilon}_t$ can be defined until the time it reaches the boundary of $[-\sqrt{\varepsilon},\varepsilon] \times \mathbb{R}^{d-1}$. Let  ${\widetilde{Z}}^{z, \varepsilon}_t $ be the first coordinate of ${\widetilde{X}}^{z, \varepsilon}_t$ and  ${{\widehat{Z}}}^{z, \varepsilon}_t $ be the vector consisting of the remaining $d-1$ coordinates of ${\widetilde{X}}^{z, \varepsilon}_t$. Note that the process $({\widetilde{Z}}^{z, \varepsilon}_t, {{\widehat{Z}}}^{z, \varepsilon}_t )$ can be written in the form (\ref{mij})-(\ref{mijxx}) with the coefficients satisfying (\ref{ez}), (\ref{mij2}), provided that $r$ and $\varepsilon$ chosen to be sufficiently small, independently of $x$.

By (\ref{iui}),
\begin{equation} \label{yuy}
\lim_{\varepsilon \downarrow 0} \mathrm{P} \left( \min( \tau^{x,\varepsilon}(S_{-\sqrt{\varepsilon}}), \tau^{x,\varepsilon}(S_\varepsilon)) > \tau^{x,\varepsilon}({U}^r_\varepsilon(x)) \right) = 0,
\end{equation}
and it is not difficult to see that the limit is uniform in $x \in \partial D$. Since $\alpha$ is bounded from above and below, from (\ref{us2}) and (\ref{us1}) it follows that there are  $\varepsilon_0 >0$ and $L > 0$ such that
\begin{equation} \label{us2xgzz}
\mathrm{E}\left( \lambda(t: X^{x,\varepsilon}_t \in \overline{U}, t \leq \tau^{x,\varepsilon}({U}^r_\varepsilon(x)))\right) \leq L \varepsilon^{2},~~~x \in \partial D,
\end{equation}
provided that $\varepsilon \leq \varepsilon_0$, and that
\begin{equation} \label{us1xgzz}
\lim_{\varepsilon \downarrow 0} \mathrm{P}\left( \tau^{x,\varepsilon}(S_\varepsilon) < \tau^{x,\varepsilon}(S_{-\sqrt{\varepsilon}})\right) =
(1+ 2  a(x))^{-1},~~~x \in \partial D.
\end{equation}
The convergence is uniform in $x$ since the dependence of the process $({\widetilde{Z}}^{z, \varepsilon}_t, {{\widehat{Z}}}^{z, \varepsilon}_t )$ on $x$ manifests itself through the value of $ \mathrm{v} = a(x)$ and through the values of $C$, $\varkappa(\varepsilon)$ in (\ref{ez}), (\ref{mij2}) in a way that doesn't affect the applicability of (\ref{us2}) and (\ref{us1}).

The  strong Markov property of the process $X^{x, \varepsilon}_t$, together with (\ref{interm}), (\ref{brmo}), and (\ref{yuy}), allows us to obtain (\ref{us2xg}) from (\ref{us2xgzz}) and (\ref{us1xg}) from (\ref{us1xgzz}).
\\


%

4) Let us get a bound on $ \mathrm{E} \lambda(u:  u \leq \sigma^{x, \varepsilon}_{1}, X^{x, \varepsilon}_u
\in \overline{U})$, $x \in \partial D$.
Since the process is pure Brownian motion on $U$,  there are  $\varepsilon_0 >0$ and $L > 0$ such that
\begin{equation} \label{us2xgkk}
\mathrm{E}\min( \tau^{x,\varepsilon}(\partial D), {\tau}^{x,\varepsilon}(S_{\sqrt{\varepsilon}})) \leq L \varepsilon^{\frac{3}{2}},~~~x \in S_{\varepsilon},
\end{equation}
provided that $\varepsilon \leq \varepsilon_0$, while
\begin{equation} \label{us1xgkk}
\lim_{\varepsilon \downarrow 0} ( \varepsilon^{-\frac{1}{2}} \mathrm{P}(\tau^{x,\varepsilon}(S_{\sqrt{\varepsilon}})  <  {\tau}^{x,\varepsilon}(\partial D) )) = 1,~~~{\rm uniformly}~{\rm in}~x \in S_{\varepsilon}.
\end{equation}
Let $\delta$ be sufficiently small for (\ref{us2xg}) and (\ref{us1xg}) to hold.
Formulas (\ref{us2xgkk}),   (\ref{us1xgkk}) together with    (\ref{us2xg}),  (\ref{us1xg}), and  the strong Markov property of the process, imply that there are  $\varepsilon_0 >0$ and $L > 0$ such that for all sufficiently small $\varepsilon$,
\begin{equation} \label{lok}
\mathrm{E} \lambda(u:  u \leq \sigma^{x, \varepsilon}_{1}, X^{x, \varepsilon}_u
\in \overline{U}) \leq L \varepsilon,~~~x \in \partial D.
\end{equation}
Let $\xi^{x, \varepsilon}_n = \lambda(u: \tau^{x, \varepsilon}_{n} \leq u \leq \sigma^{x, \varepsilon}_{n}, X^{x, \varepsilon}_u
\in \overline{U})$. Formulas (\ref{lok}) and (\ref{jj1}), together with the strong  Markov property of the process, imply that there is $c= c(t)$ such that
\begin{equation} \label{bound2}
 \mathrm{E}  \sum_{n = 0}^\infty \xi^{x, \varepsilon}_{n+1} \chi_{ \{\sigma^{x, \varepsilon}_{n} \leq s(t) \}} \leq c \sqrt{\varepsilon},~~~x \in \overline{U}.
\end{equation}
\\

5) We claim that for sufficiently small $\delta$ there are  $\varepsilon_0 >0$ and $L > 0$ such that
\begin{equation} \label{qq2}
\mathrm{E}\left( |X^{x, \varepsilon}_{\tau^{x,\varepsilon}(S_\varepsilon)} - x|^2,~ \tau^{x,\varepsilon}(S_\varepsilon) < \tau^{x,\varepsilon}(S_{-\delta})\right) \leq L {\varepsilon^2},~~x \in \partial D,
\end{equation}
provided that $\varepsilon \leq \varepsilon_0$.  Let us sketch a proof of this statement. First, by observing the process in the $\varepsilon$-neighborhood of $\partial D$, it is easy to see that
\begin{equation} \label{qq2xz}
\mathrm{E} |X^{x, \varepsilon}_{\tau^{x,\varepsilon}(S_\varepsilon) \wedge \tau^{x,\varepsilon}(S_{-\varepsilon})} - x|^2 \leq L {\varepsilon^2},~~x \in \partial D.
\end{equation}
From the the presence of the strong drift inside $D$, it follows that there is $c > 0$ such that
\begin{equation} \label{qq2ww}
\mathrm{P}\left( \tau^{x,\varepsilon}(\partial D) > t \varepsilon^2,~ \tau^{x,\varepsilon}(\partial D) < \tau^{x,\varepsilon}(S_{-\delta})\right) \leq e^{-c t},~~x \in \partial S_{-\varepsilon},~t \geq 0,
\end{equation}
and, consequently,
\begin{equation} \label{qq2mm}
\mathrm{E}\left( |X^{x, \varepsilon}_{\tau^{x,\varepsilon}(\partial D)} - x|^2,~ \tau^{x,\varepsilon}(\partial D) < \tau^{x,\varepsilon}(S_{-\delta})\right) \leq L {\varepsilon^2},~~x \in \partial S_{-\varepsilon}.
\end{equation}
Also note that there is $c > 0$ such that
\begin{equation} \label{qq2wwd}
\mathrm{P}\left(\tau^{x,\varepsilon}(\partial D) < \tau^{x,\varepsilon}(S_{-\delta})\right) \leq 1- c,~~x \in \partial S_{-\varepsilon}.
\end{equation}
By considering consecutive visits of the process to $S_{-\varepsilon}$ and $\partial D$, employing (\ref{qq2xz}), (\ref{qq2mm}), (\ref{qq2wwd}), and using the strong Markov property, we obtain (\ref{qq2}).

Take the compact set $K \subset D$ such that $\partial K = S_{-\delta}$, where $\delta$ is sufficiently small for (\ref{us2xg}), (\ref{us1xg}), and (\ref{qq2}) to hold. Observe that, since the process is pure Brownian motion in $U$, by (\ref{us2xgkk}),
\begin{equation} \label{qq1}
\mathrm{E}|X^{x, \varepsilon}_{\tau^{x, \varepsilon}(\partial D) \wedge \tau^{x, \varepsilon}(S_{\sqrt{\varepsilon}})} - x|^2   = d \cdot \mathrm{E}({\tau^{x, \varepsilon}(\partial D) \wedge \tau^{x, \varepsilon}(S_{\sqrt{\varepsilon}})}) \leq L {\varepsilon^{\frac{3}{2}}} ,~~x \in S_\varepsilon.
\end{equation}
Therefore, by (\ref{qq2}), (\ref{us1xg}), and (\ref{us1xgkk}), it follows from the strong  Markov property that
\begin{equation} \label{mjj}
\lim_{\varepsilon \downarrow 0}(\varepsilon^{-\frac{1}{2}} \mathrm{P} (\sigma^{x, \varepsilon}_1 < \tau^{x, \varepsilon}(K))) = \sum_{n = 1}^\infty (\frac{1}{1 + 2 a(x)})^n  = \frac{1}{2 a(x)},~~~{\rm uniformly}~{\rm in}~x \in \partial D.
\end{equation}
\\

6) 
Combining  (\ref{us1xg}), (\ref{us1xgkk}),  (\ref{qq2}), and (\ref{qq1}), and using the strong  Markov property, we obtain that there are  $\varepsilon_0 >0$ and $L > 0$ such that
\begin{equation} \label{nni}
 \mathrm{E} (|X^{x, \varepsilon}_{\sigma^{x, \varepsilon}_1} -x|^2, \sigma^{x, \varepsilon}_1 < \tau^{x, \varepsilon}(K)) \leq L {\varepsilon^{\frac{3}{2}}},~~x \in \partial D,
\end{equation}
provided that $\varepsilon \leq \varepsilon_0$. Therefore,
\begin{equation} \label{nnizz}
 \mathrm{E} (|X^{x, \varepsilon}_{\sigma^{x, \varepsilon}_1} -x|, \sigma^{x, \varepsilon}_1 < \tau^{x, \varepsilon}(K)) \leq L {\varepsilon},~~x \in \partial D
\end{equation}
provided that $\varepsilon \leq \varepsilon_0$, with a different constant $L$.
\\

7) Let $\overline{f}$ be the value of $f$ on $\partial D$.  Let us show that
\begin{equation} \label{jj5}
\lim_{\varepsilon \downarrow 0} (\varepsilon^{-\frac{1}{2}}  \mathrm{E} (f(X^{x, \varepsilon}_{\sigma^{x, \varepsilon}_1}) - \overline{f}))  = 0,~~~{\rm uniformly}~{\rm in}~x \in \partial D.
\end{equation}
Introduce the following two sequences of stopping times:   $\overline{\tau}^{x, \varepsilon}_1 = \tau^{x, \varepsilon}(\partial D)$, $\overline{\sigma}^{x, \varepsilon}_n = \inf(t \geq \overline{\tau}^{x, \varepsilon}_n: X^{x, \varepsilon}_t \in K)$, $n \geq 1$, while $\overline{\tau}^{x, \varepsilon}_n = \inf(t \geq \overline{\sigma}^{x, \varepsilon}_{n-1}: X^{x, \varepsilon}_t \in \partial D)$, $n \geq 2$. Then
\[
 \mathrm{E} f(X^{x, \varepsilon}_{\sigma^{x, \varepsilon}_1}) - \overline{f} =  \sum_{n =1}^\infty \mathrm{E} ( (f(X^{x, \varepsilon}_{\sigma^{x, \varepsilon}_1}) -\overline{f}), \overline{\tau}^{x, \varepsilon}_n < \sigma^{x, \varepsilon}_1 < \overline{\tau}^{x, \varepsilon}_{n+1})=
\]
\[
\mathrm{E} ( (f(X^{x, \varepsilon}_{\sigma^{x, \varepsilon}_1}) -\overline{f}), \sigma^{x, \varepsilon}_1 < \overline{\tau}^{x, \varepsilon}_{2}) +
\sum_{n =2}^\infty \mathrm{E} ( (f(X^{x, \varepsilon}_{\sigma^{x, \varepsilon}_1}) -\overline{f}), \sigma^{x, \varepsilon}_1 < \overline{\tau}^{x, \varepsilon}_{n+1} | \sigma^{x, \varepsilon}_1 > \overline{\sigma}^{x, \varepsilon}_{n-1} ) \mathrm{P}( \sigma^{x, \varepsilon}_1 > \overline{\sigma}^{x, \varepsilon}_{n-1} ).
\]
By (\ref{mjj}) and the strong  Markov property of the process, there is $c > 0$ such that
\[
\mathrm{P}( \sigma^{x, \varepsilon}_1 > \overline{\sigma}^{x, \varepsilon}_{n-1} ) \leq (1 - c\sqrt{\varepsilon})^{n-1}.
\]
Also observe that
\[
|\mathrm{E} ( (f(X^{x, \varepsilon}_{\sigma^{x, \varepsilon}_1}) -\overline{f}), \sigma^{x, \varepsilon}_1 < \overline{\tau}^{x, \varepsilon}_{n+1} | \sigma^{x, \varepsilon}_1 > \overline{\sigma}^{x, \varepsilon}_{n-1} )| \leq \sup_{x \in K} |\mathrm{E} ( (f(X^{x, \varepsilon}_{\sigma^{x, \varepsilon}_1}) -\overline{f}), \sigma^{x, \varepsilon}_1 < \overline{\tau}^{x, \varepsilon}_{2}  )| \leq
\]
\[
\sup_{x \in  K} |\mathrm{E} ( \langle \nabla f(X^{x, \varepsilon}_{\overline{\tau}^{x, \varepsilon}_1}),  X^{x, \varepsilon}_{\sigma^{x, \varepsilon}_1} -X^{x, \varepsilon}_{\overline{\tau}^{x, \varepsilon}_1} \rangle, \sigma^{x, \varepsilon}_1 < \overline{\tau}^{x, \varepsilon}_{2} ) | + C \sup_{x \in \partial D} \mathrm{E} (   |X^{x, \varepsilon}_{\sigma^{x, \varepsilon}_1} -x|^2, \sigma^{x, \varepsilon}_1 < \overline{\tau}^{x, \varepsilon}_{2}),
\]
where $C$ depends on the $C^2(\overline{U})$-norm of $f$.
The second term on the right hand side is bounded by $C {\varepsilon^{\frac{3}{2}}}$, with a different constant $C$, using (\ref{nni}). In order to estimate the first term, we notice that if $x \in \partial D$ and $y \in S_{\sqrt{\varepsilon}}$, then, since $\nabla f(x)$ is orthogonal to $\partial D$,
\[
| \langle \nabla f(x), y-x \rangle + \sqrt{\varepsilon} \langle \nabla f(x), n(x) \rangle| \leq c |x - y|^2
 \]
for some $c > 0$, where $n(x)$ is the unit inward (with respect to $D$) normal  to the boundary at $x$. Therefore, for $x \in K$,
\[
|\mathrm{E} ( \langle \nabla f(X^{x, \varepsilon}_{\overline{\tau}^{x, \varepsilon}_1}),  X^{x, \varepsilon}_{\sigma^{x, \varepsilon}_1} -X^{x, \varepsilon}_{\overline{\tau}^{x, \varepsilon}_1} \rangle, \sigma^{x, \varepsilon}_1 < \overline{\tau}^{x, \varepsilon}_{2} ) + \sqrt{\varepsilon} \mathrm{E} ( \langle \nabla f(X^{x, \varepsilon}_{\overline{\tau}^{x, \varepsilon}_1}), n(X^{x, \varepsilon}_{\overline{\tau}^{x, \varepsilon}_1}) \rangle, \sigma^{x, \varepsilon}_1 < \overline{\tau}^{x, \varepsilon}_{2} )| \leq
\]
\[
C \mathrm{E} ( |  X^{x, \varepsilon}_{\sigma^{x, \varepsilon}_1} -X^{x, \varepsilon}_{\overline{\tau}^{x, \varepsilon}_1}|^2, \sigma^{x, \varepsilon}_1 < \overline{\tau}^{x, \varepsilon}_{2} ),
\]
where $C$ depends on the $C^1(\overline{U})$-norm of $f$.
The right hand side is bounded by $C {\varepsilon^{\frac{3}{2}}}$, with a different constant $C$, using (\ref{nni}) and the strong  Markov property. Observe
that
\[
\lim_{\varepsilon \downarrow 0} \sup_{x \in  K} | \varepsilon^{-\frac{1}{2}} \mathrm{E} ( \langle \nabla f(X^{x, \varepsilon}_{\overline{\tau}^{x, \varepsilon}_1}), n(X^{x, \varepsilon}_{\overline{\tau}^{x, \varepsilon}_1}) \rangle, \sigma^{x, \varepsilon}_1 < \overline{\tau}^{x, \varepsilon}_{2} ) | = 0
\]
by (\ref{convqq}), (\ref{ortt}), (\ref{mjj}), and the strong Markov property of the process.  By (\ref{nnizz}), we also have
\[
\lim_{\varepsilon \downarrow 0} \sup_{x \in  \partial D} | \varepsilon^{-\frac{1}{2}}
\mathrm{E} ( (f(X^{x, \varepsilon}_{\sigma^{x, \varepsilon}_1}) -\overline{f}), \sigma^{x, \varepsilon}_1 < \overline{\tau}^{x, \varepsilon}_{2})| = 0.
\]
Combining the estimates above, we obtain (\ref{jj5}).
\\

Let us examine the second term in the right hand side of (\ref{vv1}). From (\ref{bound2}) and the boundedness of $\Delta f$ it follows that
\[
\lim_{\varepsilon \downarrow 0} \mathrm{E} \sum_{n = 1}^\infty  \int_{\tau^{x, \varepsilon}_{n} \wedge s(t)}^{\sigma^{x, \varepsilon}_n \wedge s(t)} \Delta f (X^{x, \varepsilon}_u) du =0.
\]
It remains to show that
\begin{equation} \label{rems}
\lim_{\varepsilon \downarrow 0} \mathrm{E} \sum_{n = 1}^\infty  \left(f(X^{x, \varepsilon}_{\sigma^{x, \varepsilon}_n \wedge s(t)} ) - f(X^{x, \varepsilon}_{\tau^{x, \varepsilon}_{n} \wedge s(t)}) \right) = 0.
\end{equation}
Introduce the stopping time
\[
s'(t) = \Big \{
  \begin{tabular}{ccc}
  $\sigma^{x, \varepsilon}_n~~{\rm if}~~\tau^{x, \varepsilon}_{n} < s(t) \leq \sigma^{x, \varepsilon}_n$ \\
  $s(t)~~~~~~~~~~~~~$ otherwise.
  \end{tabular}
\]
Now (\ref{rems}) will follow if we show that
\begin{equation} \label{rems2}
\lim_{\varepsilon \downarrow 0} \mathrm{E} \sum_{n = 1}^\infty  \left(f(X^{x, \varepsilon}_{\sigma^{x, \varepsilon}_n \wedge s'(t)} ) - \overline{f} \right) = 0,
\end{equation}
since the difference between (\ref{rems2}) and (\ref{rems}) is estimated from above by $2\sup_{x \in S_{\sqrt{\varepsilon}}}|f(x) - \overline{f}|$,
 which goes to zero as $\varepsilon \downarrow 0$. Let $N^{x,\varepsilon} = \max(n: \sigma^{x, \varepsilon}_n \leq s'(t)) $.  By the strong Markov property,
\[
\sup_{x \in \overline{U}}  \mathrm{E} \sum_{n = 1}^\infty  \left(f(X^{x, \varepsilon}_{\sigma^{x, \varepsilon}_n \wedge s'(t)} ) - \overline{f} \right) \leq
\sup_{x \in \overline{U}} \mathrm{E} N^{x,\varepsilon}
\sup_{x \in \partial D} \mathrm{E} (f(X^{x, \varepsilon}_{\sigma^{x, \varepsilon}_1}) - \overline{f}) .
\]
The right hand side tends to zero by (\ref{jj1}) and (\ref{jj5}).
\qed
\\
\\
{\it Proof of Theorem~\ref{mt1}.} Recall that $\Psi$ is the set of functions $\psi \in C(U')$ that have limits of all the first order derivatives as $y \in U, y \rightarrow x$, at all points $x \in \partial U$. This is a measure-defining class of functions on $U'$, i.e., if $\mu_1$ and $\mu_2$ satisfy $\int_{U'} \psi d \mu_1 = \int_{U'}\psi d\mu_2$ for every $\psi \in \Psi$, then $\mu_1 = \mu_2$.  As shown in Section~\ref{dttp}, for every $\psi \in \Psi$ and every $\lambda > 0$, there is $f \in \mathcal{D}(A)$ that satisfies $\lambda f - Af = \psi$.  We have demonstrated that (\ref{maeq}) holds for $f \in \mathcal{D}(A)$. By Lemma 3.1. in Chapter 8 of \cite{FW}, this is sufficient to guarantee the convergence if, in addition, the family  $\{Y^{x, \varepsilon}_t\}, \varepsilon > 0, x \in {U'}$ is tight. The tightness, however, is clear since the processes coincide with Brownian motion inside $U$, while all the points of $\partial U$ are identified.
\qed
\section{Applications, generalizations, and remarks} \label{appli}
\subsection{The behavior of the process at exponential time scales} \label{tone}

Let us now discuss the behavior, as $\varepsilon \downarrow 0$, of the original process $X^{x, \varepsilon}_t $ (rather than its trace $Y^{x, \varepsilon}_t $ on $\overline{U}$). If the process starts in a small neighborhood of $O_k$, then it takes time of order $\exp(V_k/\varepsilon)$ (in the sense of logarithmic equivalence) for it to reach $\partial D_k$, where $V_k = \inf_{x \in \partial D_k} V_k(x)$ and $V_k(x)$ is the quasi-potential defined in (\ref{quap}). Thus it is reasonable to study the behavior of $X^{x, \varepsilon}_t $ at exponential time scales, i.e., at times of order $\exp(\lambda/\varepsilon)$ with fixed $\lambda$.

The transitions between small neighborhoods of the equilibriums are governed by the matrix $V_{ij}$ defined in (\ref{vij}). In our case, $V_{ij} = V_i$ for all $i, j$, as explained in the Introduction. Because of this ``rough symmetry" (\cite{FF}), the notion of a metastable state (see (\ref{mesta})) should be replaced by that of a metastable distribution between the equilibriums.

To describe the metastable distribution for a given initial point $x \in \mathbb{T}^d$ and time scale $\exp(\lambda/\varepsilon)$ with $\lambda > 0$, assume that $V_1 < V_2 < ... < V_n$, and put $V_0 = 0$, $V_{n+1} = \infty$. We introduce the following non-standard boundary problem, which will be referred to as the $(k,j)$-problem on $U$. Namely, for each $1 \leq k \leq n$ and $k \leq j \leq n$, let
$u_{k,j}$ solve
\[
\Delta u_{k,j}(x) = 0,~~x \in U~;
\]
\[
u_{k,j}(x) = c^i_{k,j},~~x \in \partial D_i~~~{\rm and}~~~ \int_{\partial D_i} \langle \nabla u_{k,j}(x) , n(x) \rangle \nu_i(d x) =0~~~~{\rm for}~1 \leq i < k~;
\]
\[
u_{k,j}(x) = 0,~~x \in \partial D_i,~~~~{\rm for}~~i \geq k, i \neq j~;
\]
\[
u_{k,j}(x) = 1,~~x \in \partial D_j~.
\]
The constants $c^i_{k,j}$ are not prescribed, i.e., solving the $(k,j)$-problem includes finding the boundary values of the function on $\partial D_i$, $i < k$. As in Section~\ref{dttp}, the solution exists and is unique in $C^2( \overline{U})$.
\begin{theorem} \label{mstb}
Assume that $V_{k-1} < \lambda < V_k$ for some $1 \leq k \leq n$. Suppose that $T(\varepsilon)$ is such that $\lim_{\varepsilon \downarrow 0} (\varepsilon \ln T(\varepsilon))  = \lambda $. Let $\mathcal{E}_i$, $i =1,...,n$, be arbitrary disjoint neighborhoods of~$O_i$, $i =1,...,n$. Then
\[
\lim_{\varepsilon \downarrow 0} \mathrm{P}(X^{x, \varepsilon}_{T(\varepsilon)} \in \mathcal{E}_j )= u_{kj}(x),~~~x \in U,~~j \geq k~;
\]
\[
\lim_{\varepsilon \downarrow 0} \mathrm{P}(X^{x, \varepsilon}_{T(\varepsilon)} \in \mathcal{E}_j )= c^i_{k,j},~~~x \in \overline{D}_i,~~i < k,~j \geq k~;
\]
\[
\lim_{\varepsilon \downarrow 0} \mathrm{P}(X^{x, \varepsilon}_{T(\varepsilon)} \in \mathcal{E}_j ) = 1,~~~x \in \overline{D}_j,~~j \geq k,
\]
where $u_{k,j}$ is the solution of the $(k,j)$-problem and $c^i_{k,j}$ is the value of $u_{k,j}$ on $\partial D_i$. If $\lambda \geq V_n$, then
\[
\lim_{\varepsilon \downarrow 0} \mathrm{P}(X^{x, \varepsilon}_{T(\varepsilon)} \in \mathcal{E}_n ) = 1,~~x \in \mathbb{T}^d.
\]
\end{theorem}

Observe that $\sum_{j = k}^n u_{kj}(x) = 1$ for $x \in U$, while $\sum_{j = k}^n c^i_{k,j} = 1$ for $x \in \overline{D}_i$ if $i < k$. Thus, with probability close to one, $X^{x, \varepsilon}_{T(\varepsilon)}$ is located in a small neighborhood of one of the equilibrium points.

The proof of this theorem can be derived from the following facts by using the strong Markov property of the process.

1) For each $i < k$ and $x \in D_i$, the time it takes  $X^{x, \varepsilon}_t$ to exit $D_i$ is significantly smaller than $T(\varepsilon)$, i.e.,
$\tau^{x, \varepsilon}(\partial D_i)/T(\varepsilon) \rightarrow 0$ in probability as $\varepsilon \downarrow 0$.

2) For each $i \geq k$ and $x \in D_i$, the time it takes  $X^{x, \varepsilon}_t$ to exit $D_i$ is significantly larger than $T(\varepsilon)$, i.e.,
$\tau^{x, \varepsilon}(\partial D_i)/T(\varepsilon) \rightarrow \infty$ in probability as $\varepsilon \downarrow 0$.

3) The trace of the process $X^{x, \varepsilon}_t$ in $\overline{U}$ converges to $X^x_t$ as $\varepsilon \downarrow 0$.
\\
We omit the details of the proof.

\subsection{Trapping regions with multiple equilibriums}

\begin{figure}[htbp]
%
%
\centerline{\includegraphics[height=4.2in, width= 5.2in,angle=0]{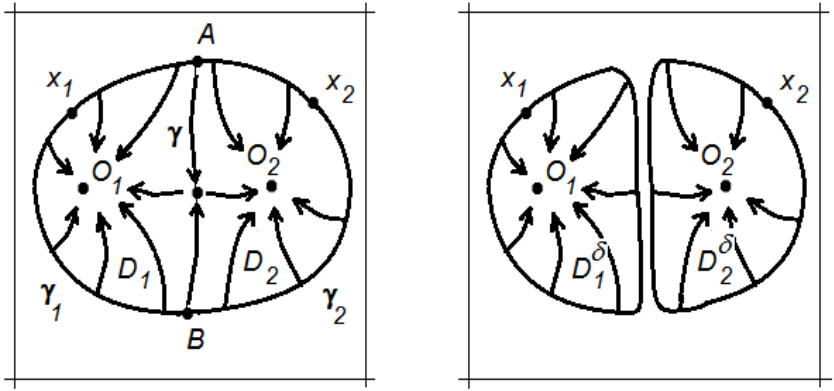}}
    \vskip -120pt
  \caption{The flow lines of the original and the modified vector fields.}
    \label{twofields}
\end{figure}

Up to now, we assumed that there was just one attractor (asymptotically stable equilibrium) inside each trapping region. Let us now consider
an example where this is not the case. For simplicity, assume that there is one trapping region $D$ containing two equilibriums $O_1$ and $O_2$
and one saddle point $S$. The structure of the vector $v$ field on $\overline{D}$ is assumed to be as shown in Figure~\ref{twofields}. As before, the vector field is equal to zero in $U = \mathbb{T}^2 \setminus \overline{D}$.

Let $D_k\subset D$ and $\gamma_k \subset \partial D$, $k =1,2$, be the sets of points that are carried to an arbitrarily small neighborhood of $O_k$ by the
deterministic flow $\dot{x}(t) = v(x(t))$. Let $A, B \in \partial D$ be the points that separate $\gamma_1$ from $\gamma_2$. Let $\gamma$ be the curve that connects $A$ with $B$ and consists of two flow lines and the saddle point (see Figure~\ref{twofields}). The asymptotic behavior  of the process $X^{x, \varepsilon}_t $ (in exponential time scales) and of the trace $Y^{x, \varepsilon}_t $ is determined by the numbers $V_{ij}$ defined in (\ref{vij}) and by the values $V_k = \inf_{x \in \partial D_k} V_k(x)$, where the quasi-potentials $V_k(x)$ are defined in (\ref{quap}).

  \vskip -10pt
  \vskip 10pt

Consider the case when $V_1 < V_{12}$, $V_2 < V_{21}$, and the infimum in the definition of $V_k$ is achieved in a unique point $x_k \in \gamma_k$, $k = 1,2$.  Then, with probability that tends to one as $\varepsilon \downarrow 0$, the process exits $D$ in an arbitrarily small neighborhood of $x_k$ provided that it starts in a small neighborhood of $O_k$, $k=1,2$.

For $\delta > 0$, one can consider the following auxiliary system. Let $V_\delta$ be the $\delta$-neighborhood of $\gamma$. Let $D^\delta_k \subset D_k$, $k =1,2$, be domains with smooth boundaries such that $D^\delta_k \setminus V_\delta =  D_k \setminus V_\delta$, yet $\overline{D^\delta_k} \bigcap \gamma = \emptyset$. Moreover, we can modify the vector field $v$ (i.e., replace it by a new vector field ${v^\delta}$) in such a way that ${v^\delta}(x) = v(x)$ for $x \notin V_\delta$, while
the field ${v^\delta}$ satisfies the assumptions with respect to the domains $D^\delta_1$ and $D^\delta_2$ that were imposed on $v$ in Section~\ref{seconv}, i.e., it is equal to zero outside $\overline{D^\delta_1} \bigcup \overline{D^\delta_2}$, is directed inside each of the domains on the boundary, and all the points of $\overline{D^\delta_k}$ are attracted to $O_k$.

The analysis of Section~\ref{seconv} applies to the process
$X^{\delta,x, \varepsilon}_t$ defined via
\[
d X^{\delta, x, \varepsilon}_t = \frac{1}{\varepsilon} v^\delta ( X^{\delta, x, \varepsilon}_t) d t + d W_t,~~~X^{\delta, x, \varepsilon}_0 = x.
\]
Let $X^{\delta, x}_t$ denote the limit of the trace process in $\mathbb{T}^2 \setminus (D^\delta_1 \bigcup D^\delta_2)$. Since $x_1, x_2 \notin V_\delta$ for all sufficiently small $\delta$, it is not difficult to show that  the probability that $X^{\delta, x}_t$ enters $V_\delta$ prior to time $T$ tends to zero for each finite $T$. Moreover, there exists the limit in probability as $\delta \downarrow 0$ of $X^{\delta, x}_t$. This limit will be denoted by $X^x_t$. This process is the limit, as $\varepsilon \downarrow 0$ of the trace $Y^{x, \varepsilon}_t$ of the original process $X^{x, \varepsilon}_t$. A direct construction of the process $X^x_t$ (in terms of the generator rather that via approximating processes) seems to be technically complicated and is not presented here.

Another case where the limit of the trace process can be easily described is when $V_1 < V_{12}$ if we assume that the infimum in the definition of $V_1$ is achieved in a unique point $x_1 \in \gamma_1$, while the infimum in the definition of $V_2$ is achieved in a unique point $x_2 \in \gamma \setminus \partial D$ (which implies that $ V_2 = V_{21}$). Thus, with probability that tends to one as $\varepsilon \downarrow 0$, the process exits $D$ in an arbitrarily small neighborhood of $x_1$ irrespective of whether it starts in $D_1$ or $D_2$. The results of Section~\ref{seconv} then apply, with the limit of the exit measure $\mu$ being the point mass concentrated at $x_1$.

A more general situation of several equilibriums  within $D$ with various relations on the quantities $V_{ij}$ and $V_k$ can be analyzed using the construction above based on removing the $\delta$-neighborhoods of the boundaries of $D_k$ and the results of Sections 6.5-6.6 of \cite{FW} on the hierarchies of cycles.

\subsection{Other generalizations}
If the process $ X^{x, \varepsilon}_t$ is governed by a more general elliptic operator
\begin{equation} \label{ope1}
L^\varepsilon = \frac{1}{2} \sum_{i,j =1}^d a_{ij}(x) \frac{\partial^2}{\partial x_i \partial x_j} + \sum_{i =1}^d b_i(x) \frac{\partial}{\partial x_i} + \frac{1}{\varepsilon} v(x) \nabla,
\end{equation}
then the results and the proofs are similar. The definition of the numbers $V_{ij}$ and of the quasi-potential $V_k(x)$ should now be based on the action functional corresponding to the operator $L^\varepsilon$. The definition (\ref{meas1}) of the measures $\nu_k$ needs to be modified to account for the variable diffusion coefficients of the process $ X^{x, \varepsilon}_t$. However, if the infimum of  $V_k(x), x \in \partial D_k$,
is achieved in a single point $x_k$, then $\nu_k$ is still the $\delta$-measure concentrated at $x_k$. 

The assumptions on the vector field $v$ that we made in Section~\ref{seconv} do not specify that  $D$ is necessarily contains a single equilibrium point. They may hold, for example, if $D$ contains a single limit cycle instead. The case of several limit cycles is technically not different from the case of several equilibrium points that we discussed above.

The results also apply to processes on general smooth manifolds, not only on a torus.

\subsection{More on the limiting process}

\begin{figure}[htbp]
%
%
\centerline{\includegraphics[height=2.1in, width= 2.7in,angle=0]{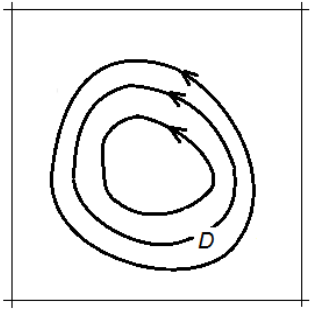}}
    \vskip -20pt
  \caption{The flow lines inside $D$.}
    \label{rot}
\end{figure}

The non-standard boundary problem introduced in Section~\ref{dttp} and the corresponding Markov process with jumps at the boundary arise in other situations, not just in the large deviation case. Consider, for example, a vector field $v$ with closed flow lines that is equal to zero outside of $\overline{D}$ such that $\partial D$ serves as one of the flow lines (see Figure~\ref{rot}).

  \vskip -20pt
  \vskip 10pt

We expect that the trace in $ \mathbb{T}^2 \setminus D$  of the process $X^{x,\varepsilon}_t$ with generator (\ref{ope1}) converges, as $\varepsilon \downarrow 0$, to the process described in Section~\ref{dttp}. The measure $\nu$ on $\partial D$ will be defined by the values of  $v$ in an arbitrarily small neighborhood of $\partial D$ and by the diffusion coefficients $a_{ij}$ on $\partial D$ and can be calculated explicitly.
\\
\\

\noindent {\bf \large Acknowledgements}: While working on this
article, M. Freidlin was supported by NSF grant DMS-1411866
and L. Koralov was supported by NSF grant DMS-1309084.
\\
\\


\begin{thebibliography}{999999}


%








\bibitem{Day} Day M., {\it Mathematical Apprach to the problem of noise induced exit}, Stochastic Analysis, Control, Optimization, and Applications.
A volume in honor of W. H. Fleming. Editors: W. McEneauey, G. Yin, Q. Zhang. Birkhauser, Boston, 1999, pp 269-287.

\bibitem{Dyn} Dynkin E. B., {\it Markov Processes}, Springer-Velag, Berlin, Heidelberg, New York, 1965.

\bibitem{FF} Freidlin M. I., {\it On Stochastic Perturbations of Systems with Rough Symmetry. Hierarchy of Markov Chains}, Journal of Statistical Physics, Vol. 157, No 6, pp 1031-1045, 2014.

\bibitem{FW} Freidlin M. I., Wentzell A. D., {\it Random
Perturbations of Dynamical Systems}, Springer 1998.

%
\bibitem{Mandl} Mandl P., {\it Analytical Treatment of One-dimensional Markov Processes}, Springer-Verlag, 1968.



\end{thebibliography}
\end{document}